\documentclass[11pt]{amsart}
\usepackage[T1]{fontenc}
\numberwithin{equation}{section}

\usepackage{float}
\usepackage{tikz-cd}
\usepackage{amsmath}
\usepackage{amsfonts}
\usepackage{amsthm}
\usepackage{amssymb}
\usepackage{graphicx}
\usepackage{hyperref}

\usepackage{float}
\usepackage{tikz-cd}
\usepackage{amsmath}%
\usepackage{amsfonts}%
\usepackage{amsthm}
\usepackage{amssymb}%
\usepackage{graphicx}
\usepackage{hyperref}
\usepackage{enumerate}
\usepackage{comment}

\usepackage{geometry}
\usepackage{layout}
\usepackage[utf8]{inputenc}
\usepackage{amsmath}
\usepackage{graphicx}
\usepackage{amsmath}
\usepackage{amssymb}
\usepackage{amsfonts}
\usepackage{graphicx}
\usepackage{url}
\usepackage{MnSymbol}
\usepackage{arydshln}
\usepackage{outlines}
\usepackage{tikz-cd}

\usepackage{subfig}
\usepackage{amsthm}

\newcommand{\Z}{\mathbb{Z}}

\newcommand{\cA}{\mathcal{A}}
\newcommand{\cO}{\mathcal{O}}
\newcommand{\cH}{\mathcal{H}}
\newcommand{\bZ}{\overline{Z}}
\newcommand{\cM}{\mathcal{M}}
\newcommand{\cS}{\mathcal{S}}
\newcommand{\cJ}{\mathcal{J}}
\newcommand{\F}{\mathbb{F}}

\renewcommand{\epsilon}{\varepsilon}

\newcommand{\cMbar}{\overline{\mathcal{M}}}
\renewcommand{\cJ}{\mathcal{J}}
\newcommand{\Jac}{\mathcal{J}}

\renewcommand{\o}{\varnothing}

\renewcommand{\o}{\varnothing}

\newcommand{\rank}{\mathrm{rank}}
\newcommand{\coker}{\mathrm{coker}}
\newcommand{\im}{\mathrm{im}}

\newcommand{\Hom}{\mathrm{Hom}}

\newcommand{\prank}{p\text{-}\rank}

\usepackage{multirow}

\usepackage{xcolor}

\makeatletter
\newtheorem*{rep@thm}{\rep@title}
\newcommand{\newrepthm}[2]{%
\newenvironment{rep#1}[1]{%
 \def\rep@title{#2 \ref{##1}}%
 \begin{rep@thm}}%
 {\end{rep@thm}}}
\makeatother

\newtheorem{thmx}{Theorem}

\newtheorem{thm}{Theorem}[section]

\newtheorem{lem}[thm]{Lemma}
\newtheorem{prop}[thm]{Proposition}
\newtheorem{cor}[thm]{Corollary}

\newrepthm{thm}{Theorem}

\theoremstyle{definition}
\newtheorem{dfn}[thm]{Definition}

\theoremstyle{remark}
\newtheorem{rem}[thm]{Remark}
\newtheorem{exmp}[thm]{Example}

\begin{document}

\author{Du\v san Dragutinovi\'c}
\address{Mathematical Institute, Utrecht University, 3508 TA Utrecht, The Netherlands}
\email{d.dragutinovic@uu.nl}

\title[Ekedahl-Oort types of stable curves]{Ekedahl-Oort types of stable curves}

\maketitle

\begin{abstract}
We discuss an intrinsic definition of Ekedahl-Oort types for stable curves via Moonen's Hasse-Witt triples, previously used for smooth curves. We show that this leads to an equivalent definition to the classical one for the Ekedahl-Oort types of generalized Jacobians as semi-abelian varieties. This approach enables us to compute bounds on the dimensions of many new Ekedahl-Oort loci of curves and to generalize results on the dimensions of the $\prank$ and $a$-number loci.
\end{abstract}

\section{Introduction}
The purpose of this article is to discuss the interpretation of the Ekedahl-Oort types of stable curves via Moonen's Hasse-Witt triples (\cite{moonen}), show the equivalence with the classical notions, and use this instrinsic description to present the full potential of the inductive technique from \cite[Theorem 2.3]{fabervdgeer}. Namely, we compute the bounds on the dimensions of Ekedahl-Oort loci satisfying certain technical conditions, and in particular, of all ${\prank \leq f}$ and $a$-number $\geq a$ loci of $\cMbar_g$, the moduli space of stable curves of genus $g$ in characteristic $p>0$.
This technique has been used in \cite{fabervdgeer}, \cite{achterpries}, and \cite{pries_a_number} to describe and compute the dimensions of the components of the $\prank \leq f$ loci in $\cMbar_g$ and to bound the dimensions of some $a$-number $\geq 2$ subloci, and similarly in \cite{glasspries} and \cite{achterprieshe} for hyperelliptic curves. Let us introduce some notation and briefly present the results.

 Let $k$ be an algebraically closed field of characteristic $p>0$ and $g\geq 2$.
For a principally polarized $g$-dimensional abelian variety $A$ over $k$, the $\prank$ and the {$a$-number} of $A$, are the numbers $f(A)$ and $a(A)$, defined by 
\begin{equation}
f(A) = \dim_{\F_p}\Hom(\mu_p, A[p])\quad \text{and}\quad a(A) = \dim_k\Hom(\alpha_p, A[p]),
\label{eqn:prankanum}
\end{equation}
where $\mu_p$ and $\alpha_p$ are the kernels of the Frobenius morphism on the multiplicative group~$\mathbb{G}_m$ and the additive group $\mathbb{G}_a$, respectively, and where $\Hom(\alpha_p, A[p])$ is viewed as a $k$-vector space via the identification $k \cong \mathrm{End}_k(\alpha_p)$; alternatively, $f(A) = f$ if $\# A[p](k) = p^f$. The isomorphism class of the $p$-torsion group scheme $A[p]$ corresponds to the Young diagram  $$\mu = [\mu_1, \mu_2, \ldots, \mu_n], $$ for some $g\geq \mu_1 > \ldots > \mu_n > 0$; we describe this in Section \ref{sec:eo_types}. We call $\mu$ the \textit{Ekedahl-Oort type} of $A$, and write $\mu(A) = \mu$. Note that $f(A) = g - \mu_1$  and $a(A) = n$. These definitions extend to the case where $A$ is a semi-abelian variety. 

Let $\cM_g$ denote the moduli space of smooth genus-$g$ curves and $\cMbar_g$ the moduli space of stable genus-$g$ curves, its Deligne-Mumford compactification.
The boundary $\cMbar_g - \cM_g$ consists of the components $\Delta_0, \Delta_1, \ldots, \Delta_{\left \lfloor \frac{g}{2} \right \rfloor}$. The generic point of $\Delta_0$ corresponds to an irreducible curve $C_0$ with an ordinary double point, whose normalization $\Tilde{C}_0$ is a genus-$(g - 1)$ curve, while the generic point of $\Delta_i$ corresponds to a reducible curve with components of genus $i$ and $g - i$ for any $1 \leq i \leq \left \lfloor \frac{g}{2} \right \rfloor$.
To a smooth curve $C$ of genus $g$, we can attach its Jacobian $\Jac_C = \mathrm{Pic}^0(C)$. This induces the Torelli morphism $j: \cM_g \to \cA_g$, where  $\cA_g$ is the moduli space of principally polarized $g$-dimensional abelian varieties. For a stable curve $C$, its generalized Jacobian $\Jac_C = \mathrm{Pic}^0(C)$ is a semi-abelian variety. By \cite[Theorem 4.1]{alexeev}, this extends the Torelli morphism to 
\begin{equation}
 \label{eqn:torelli}
 j: \cMbar_g \to \Tilde{\cA}_g,
\end{equation}
where $\Tilde{\cA}_g$ is a certain fixed smooth toroidal compactification of $\cA_g$.

Although the invariants of curves $C$ in characteristic $p > 0$ are often defined in terms of their Jacobians $\mathcal{J}_C$, intrinsic descriptions can provide valuable insights. Alternative definitions of the $p$-rank and $a$-number of $C$ are $f(C) = \rank_k(F_C^g)$ and $a(C) = g - \rank_k(F_C)$, where $F_C$ denotes the Frobenius operator on $H^1(C, \mathcal{O}_C)$. The Ekedahl-Oort type of $C$ can also be defined in terms of its first de Rham cohomology $H^1_{dR}(C)$. In~\cite{moonen}, Moonen introduces the concept of a \emph{Hasse-Witt triple}, which can be used to equivalently define the Ekedahl-Oort type $\mu(C)$ of a smooth curve $C$. We follow this idea and extend the interpretation to all stable curves, providing an intrinsic definition equivalent to the classical one. Moreover, this description reveals certain relations more clearly.
\begin{thmx}[Theorem \ref{thm:muCnorm}]
\label{thmx:equivalence_hwt_eo}
Let $C$ be a stable curve, $\cJ_C$ its generalized Jacobian, and $h: \Tilde{C} \to C$ its normalization. Then  $\mu(\cJ_C) = \mu(C) =  \mu(\Tilde{C})$, where $\mu(\cJ_C)$ is classically defined, while $\mu(C)$ and $\mu(\Tilde{C})$ are defined in terms of Hasse-Witt triples. 
\end{thmx}
This description enables us to achieve the full potential of the inductive technique used in \cite[Theorem 2.3]{fabervdgeer}, 
and generalize some known results. We obtain an inductive formula for computing dimensions of Ekedahl-Oort strata of curves, in the case when the defining Ekedahl-Oort type satisfies certain technical conditions. For example, these technical conditions are satisfied for all loci in $\cMbar_g$ of curves with $p$-rank $\leq f$ and $a$-number~$\geq a$, where $1\leq a + f \leq g$. We denote by $\bZ_{\mu}(\cMbar_g)$ the preimage of the closed locus $\bZ_{\mu} \subseteq \Tilde{\cA}_g$ under the Torelli morphism \eqref{eqn:torelli}, and $\bZ_{\mu}$ is the closure of the locus of principally polarized abelian varieties of dimension $g$ in characteristic $p$ with Ekedahl-Oort type $\mu$. 
\begin{thmx}[Theorem \ref{thm:main1}]
\label{thmx:induction_eo_curves}
Let $g\geq 2$, let $\mu$ be a Young diagram that satisfies the condition~\eqref{eqn:technical_cond}, and let $d = 3g - 3 - \max \{\dim\Gamma':\Gamma' \text{ a component of }\bZ_{\mu}(\cMbar_g)\}.$ If $\Gamma\subseteq \bZ_{\mu}(\cMbar_{g + 1})$ is a component such that $\Gamma \cap \Delta_0 \neq \emptyset$, then the codimension of $\Gamma$ in~$\cMbar_{g + 1}$ is at least~$d$.
\end{thmx}

This theorem serves as a common generalization of \cite[Theorem 2.3]{fabervdgeer} and \cite[Proposition 3.7]{pries_a_number}. Namely, it recovers these results for the choices $\mu = [g]$ and $\mu = [g, 1]$, respectively. 
Moreover, it allows us to establish an upper bound for various loci of curves, such as the locus $ T_3(\cMbar_g) = \bZ_{[3, 2, 1]}(\cMbar_g)$  of stable genus-$g$ curves with $a$-number~$\geq 3$, a result that was previously out of reach; see Propositions \ref{prop:32_genus45} and~\ref{prop:321}.  

Beyond \cite{fabervdgeer} and \cite{pries_a_number}, this article addresses questions related to those in \cite[Theorem 6.4]{pries_current_results}, which outlines an inductive method for constructing smooth curves of genus~$g$ with a prescribed Ekedahl-Oort type $\mu$; see Remark \ref{rem:priessmootheo}. Additional references on the geometry of Ekedahl-Oort stratifications include \cite{pries_eo}, \cite{elkinpries_he}, and \cite{zhou_genus4}, the latter focusing on the locus of hyperelliptic genus-$4$ curves in characteristic $p = 3$.

Here is the outline of the paper. 
In Section \ref{sec:eo_types}, we present several equivalent definitions of the Ekedahl-Oort types for principally polarized abelian varieties and smooth curves, along with a discussion of their main properties. In Section~\ref{sec:bdry}, we prove Theorem~\ref{thmx:equivalence_hwt_eo}. This intrinsic perspective is then used in Section~\ref{sec:inductiveres} to establish Theorem~\ref{thmx:induction_eo_curves}. Furthermore, in Section \ref{sec:app}, we apply these results to obtain new findings regarding the dimensions of certain Ekedahl-Oort loci of curves, including an upper bound on the dimension of $T_3(\cMbar_g)$. Also,  in Theorem~\ref{thm:he_implies_eo}, we provide a criterion for computing the dimensions of almost all Ekedahl-Oort loci $\bZ_{\mu}(\cMbar_4)$, based on certain conclusions about the locus of hyperelliptic curves. Finally, in Section~\ref{sec:exm}, we give concrete examples that illustrate the geometry of the Ekedahl-Oort strata of curves in characteristic $p = 2$ and $p = 3$; in the latter case, our results follow from the previously mentioned criterion.

\subsection*{Acknowledgment}
The author is grateful to Carel Faber for all the discussions and valuable comments, to Valentijn Karemaker and Rachel Pries for their feedback and insightful suggestions, and to the anonymous referees for their helpful remarks. The author is supported by the Mathematical Institute of Utrecht University.

\section{Ekedahl-Oort types of abelian varieties}
\label{sec:eo_types}

Let $A$ be a $g$-dimensional principally polarized abelian variety over an algebraically closed field $k$ of characteristic $p>0$, and let $\sigma: k \to k, a \mapsto a^p$ be the Frobenius morphism on $k$.
Below, we mainly follow \cite[Section 5]{oort}, together with \cite{vdgeercycle} and \cite{moonen}.  

A \textit{polarized Dieudonn\'e module} over $k$ is a quadruple $(M, F, V, b)$ where \begin{itemize}
\item $M$ is a $2g$-dimensional $k$-vector space;
    \item $F: M \to M$ is a $\sigma$-linear operator, while $V: M \to M$ is a $\sigma^{-1}$-linear operator;
    \item $\ker(F) = \im(V)$ and $\ker(V) = \im(F)$;
    \item $b: M\times M \to k$ is a polarization, i.e., a non-degenerate alternating bilinear form such that $b(F(x), y) = b(x, V(y))^p$ for all $x, y \in M$.
\end{itemize}

Consider $\mathbb{D}_A = (\mathbb{D}_A, F, V, b) = \mathbb{D}(A[p])$, the (polarized) Dieudonné module of~$A$. 
Starting from the filtration $0 \subset \mathbb{D}_A$, we can refine it to a maximal sub-filtration 
$0 \subset N_1 \subset \ldots \subset N_g \subset N_{g + 1} \subset \ldots \subset N_{2g} = \mathbb{D}_A, $
called a final filtration, which is stable under the action of $V$ and $F^{-1}$, and such that $i = \dim N_i$; it satisfies $V(N_{2g}) = N_g$.
The \textit{final type} of $A$ is the increasing surjective map $\nu: \{0, 1, \ldots, 2g\} \to \{0, 1, \ldots, g\}$ defined by $\nu(0) = 0$, and for $1\leq i\leq 2g$:
\begin{equation}
\nu(i) = \dim_k V(N_i).  
\label{eqn:nu_def}
\end{equation}
While a final filtration of $A$ is not necessarily unique, the final type of $A$ is unique. Moreover, it is determined by the values $\nu(i)$ for $1 \leq i \leq g$, taking into account the relations $\nu(2g - i) = \nu(i) - i + g$ for $1 \leq i \leq 2g$, which follow from the properties of a final filtration. We denote this type by $\nu = (\nu(1), \ldots, \nu(g))$ and refer to it as having degree $g$.

An alternative combinatorial description uses Young diagrams $\mu = \mu(A) = [\mu_1, \mu_2, \ldots, \mu_n]$ with $g\geq \mu_1 > \ldots \mu_n > 0$ defined in terms of $\nu = \nu(A)$ by 
\begin{equation}
\mu_j = \#\{i: 1\leq i\leq g \text{ and } j\leq i - \nu(i)\}.
\label{eqn:nu_to_mu}
\end{equation}
We say that $\mu(A)$ is the \textit{Ekedahl-Oort type} of $A$; note that $f(A) = \mu_1$ and $a(A) = n$.   
The set $Z_{\mu}$ of all $A$ in $\cA_g$ such that $\mu(A) = \mu$ is called the \textit{Ekedahl-Oort statum}. 
By \cite[Corollary 11.2]{oort}, the locally closed set $Z_{\mu}$ is pure of codimension $$cd(\mu) = \sum_{i = 1}^n \mu_i$$ in~$\cA_g$. By \cite[Proposition 12.5]{oort}, if $\mu$ and $\mu'$ are two Young diagrams, then 
\begin{equation}
Z_{\mu'} \cap \bZ_{\mu} \neq \emptyset\implies Z_{\mu'} \subseteq \bZ_{\mu}  \text{ in } \cA_g.
\label{eqn:oort125}
\end{equation}
Furthermore, by \cite[Proposition 11.1]{oort}, if we introduce a partial ordering by $$\mu = [\mu_1, \ldots, \mu_n] \geq \mu'= [\mu'_1, \ldots, \mu'_m] \text{ if } n\leq m \text{ and } \mu_i \leq \mu_i' \text{ for all } 1\leq i\leq n, $$ for any $g$ as long as $g  \geq \max\{\mu_1, \mu_1'\}$, we get 
\begin{equation}
\mu' \leq \mu \quad \implies Z_{\mu'} \subseteq \bZ_{\mu} \text{ in } \cA_g.
\label{eq:finaltypes_EOclosure}    
\end{equation}
The reverse implication does not hold in general. In {Section~ \ref{sec:weyl}}, we discuss this and explain the technical conditions of Theorem~\ref{thm:main1}; see also \cite[14.3]{oort} and \cite[page 611]{ekedahlvdgeer}.

The Ekedahl-Oort stratification of $\cA_g$ for $g = 2, 3, 4$ was described in \cite[Section 4]{pries_eo}. Below, we present an example for $g = 5$. Additionally, see the beginning of Section \ref{sec:inductiveres} for a brief overview of the 'geography' of the Ekedahl-Oort loci of curves.

\begin{exmp}
\label{exmp:4132}
The locus $\bZ_{[4, 1]}\subset \cA_5$ contains $Z_{\mu}$ for any Young diagram  $\mu$ from the list $[4, 2], [4, 3]$, $[4, 2, 1]$, $[4, 3, 1], [4, 3, 2]$, $[4, 3, 2, 1], [5, 1]$, $[5, 2], \ldots$, $[5, 4, 3, 2]$, or $[5, 4, 3, 2, 1]$. On the other hand, $Z_{[4, 1]}$ is contained in the closed loci $\bZ_{\mu}$ for $\mu$ any Young diagram from $\emptyset, [1], [2], [2, 1]$, $[3], [3, 1]$, or $[4]$. Note that $Z_{[4, 1]} \not \subset \bZ_{[3, 2]}$ by \eqref{eqn:oort125} since $Z_{[3, 2]}$ and $Z_{[4, 1]}$ have the same dimension in~$\cA_5$ and are both irreducible by \cite[Theorem 11.5]{ekedahlvdgeer}.
\end{exmp}
An important notion that Moonen introduced is a \textit{Hasse-Witt triple} over a perfect field~$K$ of characteristic $p>0$, which is a triple $(Q, \Phi, \Psi)$ where \begin{itemize}
    \item $Q$ is a finite-dimensional $K$-vector space;
    \item $\Phi: Q \to Q$ is a $\sigma$-linear map;
    \item $\Psi: \ker(\Phi) \overset{\cong}{\to} \im(\Phi)^{\perp}$ is a $\sigma$-linear bijective map; 
\end{itemize} 
where $\im(\Phi)^{\perp} = \{\lambda \in Q^{\vee}: \lambda(q) = 0 \text{ for all }q \in \im(\Phi)\}$ and $\sigma: K \to K$ is the Frobenius morphism on $K$. Here, we choose $K = k$ and recall the following result.  

\begin{thm}[{\cite[Theorem 2.8; 2.5 and 2.6]{moonen}}]
\label{thm:moonen}
There is a bijection $$\left\{ \begin{matrix} 
\text{ isomorphism classes of }\\ (M, F, V, b),\text{ } \dim_k(M) = 2g \end{matrix} \right \} \longleftrightarrow \left\{ \begin{matrix} 
\text{ isomorphism classes of }\\ (Q, \Phi, \Psi),\text{ } \dim_k(Q) = g \end{matrix} \right \}, $$
given by $(M, F, V, b) \mapsto (Q, \Phi, \Psi)$ where $Q = M/\ker F$, $\Phi = \overline{F}:  Q\to Q$ is the map induced by $F$, and $\Psi(x) = b(-, F(x))$.
\end{thm}

\subsection{Ekedahl-Oort types of smooth curves} 
\label{subsec:eo_smooth}
Let $C$ and $\cJ_C$ be a smooth curve over~$k$ and its Jacobian. Classically, we say  $f(C) = f(\cJ_C)$ and $a(C) = a(\cJ_C)$, and define the Ekedahl-Oort type of $C$ as the Young diagram $\mu(C) = \mu(\cJ_C).$

Alternatively, we can define the Ekedahl-Oort type of $C$ as in \cite{moonen}. Namely, we take $Q =~H^1(C, \cO_C)$, $\Phi = F_C$ the Frobenius operator on $H^1(C, \cO_C)$, and $\Psi: \ker(F_C) \to \im(F_C)^{\perp}$ a $\sigma$-linear bijection as in Theorem \ref{thm:moonen}. Together, these form the Hasse-Witt triple $(Q, \Phi, \Psi)$, whose isomorphism class determines the Ekedahl-Oort type by Theorem~\ref{thm:moonen}.

Lastly, let us comment on an intrinsic definition of the Ekedahl-Oort type of a stable curve $C$, previously known in the literature. By \cite[Section~5]{oda}, it follows that
$$ \mathbb{D}(\cJ_C[p]) \cong H_{dR}^1(\cJ_C) \cong H_{dR}^1(C), $$
where $H_{dR}^1(C)$ and $H_{dR}^1(\cJ_C)$ denote the first de Rham cohomology of $C$ and $\cJ_C$, respectively. Therefore, one can use this identification to define the Ekedahl-Oort type of $C$ as the Young diagram associated with the quadruple $(H_{dR}^1(C), F, V, \langle \cdot, \cdot \rangle)$, computed using the procedure described above.

\subsection{On the Weyl group and the inclusion of Ekedahl-Oort loci}
\label{sec:weyl}
We conclude this section by discussing the inclusion of one Ekedahl-Oort locus in the closure of another, using an alternative combinatorial characterization of Ekedahl-Oort types. For a detailed overview of the combinatorial techniques, see \cite[Section 2]{ekedahlvdgeer}.

Let $S_g$ be the symmetric group on $g$ letters, let $W_g$ be the Weyl group defined by 
$$W_{g} = \{\omega \in S_{2g}: \omega(i) + \omega(2g + 1 - i) = 2g + 1 \text{ for }i = 1, 2, \ldots, g\} \cong S_g \ltimes (\Z/2\Z)^{g}.$$ 
The length of an element $\omega \in W_g$ is the number $$l(\omega) = \#\{i < j \leq g: \omega(i)>\omega(j)\} + \#\{i \leq j \leq g: \omega(i) + \omega(j)> 2g + 1\}.$$
For any $\omega \in W_g$, $\omega(i)$ determines $\omega(2g + 1 - i)$ for $1\leq i\leq g$, so we can write  $\omega = [a_1 a_2\ldots a_g]$ if $\omega(i) = a_i$, for $1\leq i\leq g$.
Let $\omega = [a_1a_2\ldots a_g]$ and $\omega' = [b_1b_2\ldots b_g]$. Using \cite[Lemma~2.1]{ekedahlvdgeer}, we write $\omega \leq \omega'$, which is the Bruhat-Chevalley order, if for all $1\leq i \leq d \leq g$, the $i$th largest element of $\{a_1, \ldots, a_d\}$ is $\leq$ the $i$th largest element of $\{b_1, \ldots, b_d\}$.

Let $W_I = \{\omega \in W_g: \omega\{1, 2, \ldots, g\} = \{1, 2, \ldots, g\}\} \cong S_g$ be a subgroup of $W_g$ and note that $|W_g|/|W_I| = 2^g$. If we denote by $^IW$ and $W^{I}$ the sets of minimal length coset representatives of $W_I\backslash W_g$ and $W_g/W_I$, the map $\omega \mapsto \omega^{-1}$ is an order preserving bijection of $^IW$ and $W^{I}$ by \cite[page 260]{pinkwedhornziegler} (or \cite[Lemma 1.5]{wedhorn}). There is a one-to-one correspondence between the cosets in $W_I\backslash W_g$ and the final types $\nu$ (of degree $g$). It is
induced by the map $\omega \mapsto \nu_{\omega} = (\nu_{\omega}(1), \ldots, \nu_{\omega}(g))$, with $\nu_{\omega}(i) = i - \#\{a \leq g: \omega^{-1}(a)\leq i\}$, which depends only on the cosets; see \cite[page 575]{ekedahlvdgeer} for more details. For a Young diagram $\mu$ corresponding to $\nu$ by \eqref{eqn:nu_to_mu}, let $\omega_{\mu} \in {^{I}{W}}$ be such that  $\nu_{\omega_{\mu}} = \nu$. 

If $\omega_{0, I}\in W_I$ is the element defined by $\omega_{0, I}(i) = g + 1 - i$, $1 \leq i \leq g$, by \cite[Corollary~6.5]{wedhorn}, we get the following equivalence for any two Young diagrams $\mu$ and $\mu'$ and any $p>0$:
\begin{equation}
\text{there is }u \in W_I\text{ such that } u\cdot \omega_{\mu'}\cdot (\omega_{0, I}u\omega_{0, I}) \leq \omega_{\mu} \Leftrightarrow  
Z_{\mu'} \subseteq \bZ_{\mu} \text{ in }\cA_g,  
\label{eqn:shuffle}
\end{equation}
This claim follows from \cite[Theorem 5.4]{wedhorn}; see \cite[Theorem 6.2]{pinkwedhornziegler} for the published reference. The element $u\cdot \omega_{\mu'}\cdot (\omega_{0, I}u\omega_{0, I})$ occurring in \eqref{eqn:shuffle} is called the \textit{shuffle} of $\omega_{\mu'}$.

\section{Ekedahl-Oort types of stable curves}
\label{sec:bdry}

In Section \ref{sec:eo_types}, we defined the $p$-rank, $a$-number, and Ekedahl-Oort type of a smooth curve~$C$ via its Jacobian $\cJ_C$. Although this approach applies to stable curves, there are also more intrinsic methods that provide equivalent definitions. Namely, given a stable curve $C$ of genus~$g$ over $k$, the $p$-rank and the $a$-number of $C$ can be defined as $$f(C) = \rank_k(F_C^{g}) \quad \text{and} \quad a(C) = g - \rank_k(F_C), $$ where $F_C: H^1(C, \cO_C) \to H^1(C, \cO_C)$ is the Frobenius operator on $C$. Below, we define the Ekedahl-Oort type of $C$ following the approach of \cite{moonen} for smooth curves and compare it with the classical methods.

Let $A$ be a semi-abelian variety of dimension $g$ that is an extension of an abelian variety $A_0$ by a torus of toric rank $l$. The Dieudonn\'e module of $A$, i.e., of $A[p]$ (see Section \ref{sec:eo_types}), equals 
\begin{equation}
\mathbb{D}_A = (M, F, V, b)    
\label{dfneq:evdg_semiab}
\end{equation}
where $M = M_{l, \text{ord}}\oplus M_{A_{0}}$, $F = (F_{l, \text{ord}}, F_{A_{0}})$, $V = (V_{l, \text{ord}}, V_{A_{0}})$ and $b = b_{l, \text{ord}} + b_{A_0}$, with $\mathbb{D}_{A_0} = (M_{A_0}, F_{A_0}, V_{A_0}, b_{A_0})$ and 
$\mathbb{D}_{l, \text{ord}} = (M_{l, \text{ord}}, F_{l, \text{ord}}, V_{l, \text{ord}}, b_{l, \text{ord}})\cong (k\oplus k^{\vee}, \sigma, \sigma^{-1}, b_{\text{ord}})^{\oplus l}$
the polarized Dieudonn\'e module of an ordinary $l$-dimensional principally polarized abelian variety and of $A_0$, respectively. We write $\mathbb{D}_A = \mathbb{D}_{l, \text{ord}}\oplus~\mathbb{D}_{A_0}$. In \cite[Section 6]{oort} and \cite[Section 5]{ekedahlvdgeer}, the final type $\nu_A$ of $A$ is defined as $\nu_A(i) = i$, for $1\leq i\leq l$, and $\nu_A(i + l) = \nu_{A_0}(i) + l$,~for $0\leq i\leq g - l$. Equivalently,
\begin{equation}
\mu(A) = \mu(A_0),
\label{eqn:dfn:eo_semiab_evdg}
\end{equation}
and we refer to $\mu(A)$ as the Ekedahl-Oort type of $A$.

Let now $C$ be a stable curve over $k$, $\cJ_C$ its generalized Jacobian, and $h: \Tilde{C} \to C$ the normalization of $C$. By \cite[Corollary 4.2]{oort:genjac}, $\cJ_C$ is an extension of $\cJ_{\tilde{C}}$ by a torus of toric rank $l = 1 - n + m$, where $n$ (resp. $m$) is the number of irreducible components (resp. singular points) of $C$. It follows that $\mu(\cJ_C) = \mu(\cJ_{\Tilde{C}})$.

\begin{lem}
\label{lemma:Fmodules_iso}
Let $C$ be a stable curve over $k$ with $n$ irreducible components and $m$ singular points, $l = 1 - n + m$, and let $h: \Tilde{C} \to C$ be its normalization. Then, for any $i \geq 1$, we have $\rank_k (F_C^i) = \rank_k (F_{\Tilde{C}}^i) + l$, where $F_C$ (resp. $F_{\Tilde{C}}$) is the Frobenius on $C$ (resp. $\Tilde{C}$).     
\end{lem}

\begin{proof}
Similarly as in \cite{bouw}, we get the short exact sequence of sheaves $$0 \to \cO_C \to h_{*}\cO_{\Tilde{C}} \to Q' \to 0, $$ with $Q' = h_{*}\cO_{\Tilde{C}}/\cO_C $ the quotient sheaf; note that $Q'$ is constant and supported at the nodes of $C$. The maps above are given by $$f' \mapsto f'\circ h \quad \text{and}\quad f'' \mapsto (f''(P_1) - f''(P_2))_R,$$ for  $P_1$ and $P_2$ the points over a node $R$ of $C$. Taking the long exact sequence in cohomology and using that $H^0(C, \cO_C) = k$, $H^0(\Tilde{C}, \cO_{\Tilde{C}}) = k^n$, $H^{0}(C, Q') = k^m$, and $H^1(C, Q') = 0$, we get the short exact sequence 
\begin{equation}
0 \to k^l \to H^1(C, \cO_C) \to H^1(\Tilde{C}, \cO_{\Tilde{C}}) \to 0,    
\label{eq:shortseq}
\end{equation}
where $l = 1 - n + m$ is the number of loops in $C$. This gives us an isomorphism of $k$-vector spaces $H^1(C, \cO_C) \cong k^l\oplus H^1(\Tilde{C}, \cO_{\Tilde{C}})$. However, this is not enough to get the desired conclusion. Luckily, even more is true. Consider the following diagram
\begin{equation*}
\begin{tikzcd}
0\arrow[r] & \cO_C \arrow[r] \arrow[d, "F"] &  h_{*}\cO_{\Tilde{C}} \arrow[r] \arrow[d, "F"] &  Q' \arrow[r] \arrow[d, "F"] &  0\\
0\arrow[r] & \cO_C \arrow[r]  &  h_{*}\cO_{\Tilde{C}} \arrow[r]  &  Q' \arrow[r]  &  0
\end{tikzcd}   
\end{equation*}
and let us show that all the squares in it commute. Once we consider the stalks, we see that~$F$ commutes with the integral closure for the left square since both $\cO_{C, P}$ and $(h_*\cO_{\Tilde{C}})_P$ are subrings of $k(C)$, and with the evaluation for the right one. Therefore, we obtain that 
\begin{equation}
H^1(C, \cO_C) \cong k^l\oplus H^1(\Tilde{C}, \cO_{\Tilde{C}})    
\label{eq:kFmods}
\end{equation}
is an isomorphism of $k[F]$-modules, and the iterations of the Frobenius respect \eqref{eq:shortseq}. Since the Frobenius acts as a bijection on $k$, we get $\rank_k(F_C^i) = \rank_k (F_{\Tilde{C}}^i) + l$, for any $i\geq 1$. 
\end{proof}

For $i = 1$ and $i = g$ we recover $a(C) = a(\Tilde{C})$ and $f(C) = f(\Tilde{C}) + l$, which are well-known.
Since the Frobenius $F_C$ acts on a disjoint union of smooth curves $C = \cup_i C_i$ as the collection of the Frobenius operators $(F_{C_i})_i$, if $(H^1(C_i, \cO_{C_i}), F_{C_i}, \Psi_{C_i})$ is a Hasse-Witt triple of $C_i$ obtained for any $i$ using the description from Section \ref{sec:eo_types}, 
one readily defines a Hasse-Witt triple $(H^1({C}, \cO_{{C}}), F_{{C}}, \Psi_{{C}}) = (\oplus_i H^1({C_i}, \cO_{{C_i}}), \oplus_i F_{{C_i}}, \oplus_i\Psi_{{C_i}})$ using the collection $((H^1(C_i, \cO_{C_i}), F_{C_i}, \Psi_{C_i}))_i$. (Namely, we specify a Hasse-Witt triple for every component $C_i$ leading to a single Hasse-Witt triple for $C = \cup_i C_i$.) By \cite[2.5]{moonen}, $\im(F_C)^{\perp} = \coker (F_C)^{\vee}$ for any stable curve~$C$.

\begin{dfn}
A \textit{Hasse-Witt triple of a stable curve} $C$ is a triple $(Q, \Phi, \Psi)$ over $k$ where 
\begin{itemize}
    \item $Q = H^{1}(C, \cO_C)$,
    \item $\Phi = F_C$ the Frobenius operator on $H^{1}(C, \cO_C)$, and
    \item  $\Psi = \Psi_{C}: \ker(F_C) \overset{\cong}{\to} \coker (F_C)^{\vee}$ a $\sigma$-linear bijective map such that, if ${h: \Tilde{C} \to C}$ is the normalization of $C$, the following diagram commutes
\begin{equation}
\label{diag:psi}
\begin{tikzcd}
\ker(F_C)    \arrow[rr] \arrow[d, "\Psi_{C}"]                                                                                          &  & \ker(F_{\Tilde{C}}) \arrow[ "\Psi_{\Tilde{C}}", d]                                                                                       \\
{\coker(F_C)^{\vee}} &  & {\coker(F_{\Tilde{C}})^{\vee} } \arrow[ll]
\end{tikzcd},  
\end{equation} 
where the top and the bottom maps are induced by $h$.
\end{itemize} 
\label{dfn:hw_stable}
\end{dfn}

In the proof of Theorem \ref{thm:muCnorm}, we show that the isomorphism class of a Hasse-Witt triple $(H^1(C, \cO_C), F_C, \Psi_C)$ does not depend on the choice of $\Psi_C$ as long as $\Psi_C$ is a $\sigma$-linear bijection such that \eqref{diag:psi} holds; see also \cite[2.9]{moonen}. Furthermore, note that Definition~\ref{dfn:hw_stable} is intrinsic - it does not rely on the generalized Jacobian $\cJ_C$ of $C$.

\begin{dfn} For a stable curve or a finite disjoint union of smooth curves $C$, let  $\mathbb{D}_C$ be a Dieudonn\'e module associated with the isomorphism class of any Hasse-Witt triple of $C$ by Theorem \ref{thm:moonen},  and denote by $\mu(C)$ the Young diagram associated with $\mathbb{D}_C$ as in Section~\ref{sec:eo_types}. 
We say that the \textit{Ekedahl-Oort type} of $C$ equals $\mu(C)$.
\label{dfn:eo_stable}
\end{dfn}

We have the following comparison.

\begin{thm}
Let $C$ be a stable curve over $k$ and $\cJ_C$ its generalized Jacobian. Then  $\mu(C) = \mu(\cJ_C)$, where $\mu(C)$ is as in Definition \ref{dfn:eo_stable} and $\mu(\cJ_C)$ is as in \eqref{eqn:dfn:eo_semiab_evdg}. Moreover, if $h: \Tilde{C} \to C$ is the normalization of $C$, then $\mu(C) = \mu(\Tilde{C})$, with $\mu(\Tilde{C})$ as in Definition~\ref{dfn:eo_stable}.
\label{thm:muCnorm}
\end{thm}

\begin{proof}
Let $C$ be a stable curve with $n$ irreducible components $C_i$ and $m$ singular points. Let $\Tilde{C} \to C$ be its normalization and let $l = 1 - n + m$. Lemma \ref{lemma:Fmodules_iso} describes $H^1(C, \cO_C)$ as a $k[F_C]$-module. We need to show that there is a $\sigma$-linear map $\Psi_C$ as in Definition \ref{dfn:hw_stable}, unique up to isomorphism, and that it satisfies the required properties.
 Given a Hasse-Witt triple $(H^1(\Tilde{C}, \cO_{\Tilde{C}}), F_{\Tilde{C}}, \Psi_{\Tilde{C}})$ of $\Tilde{C}$, note that $\Psi_{\Tilde{C}}$ defines $\Psi_{C}$ uniquely up to isomorphism via~\eqref{diag:psi}. Indeed, $F_C$ coincides with the bijection $\sigma$ on $k$, leading to the following isomorphism: 
$\ker(F_C) = 0^{\oplus l}\oplus \ker(F_{\Tilde{C}}) \cong \ker(F_{\Tilde{C}})$, and $${\coker(F_C) = \frac{k^{\oplus l}\oplus H^{1}(\Tilde{C}, \cO_{\Tilde{C}})}{\sigma(k)^{\oplus l}\oplus F_{\Tilde{C}}(H^{1}(\Tilde{C}, \cO_{\Tilde{C}}))} }\cong \frac{H^{1}(\Tilde{C}, \cO_{\Tilde{C}})}{F_{\Tilde{C}}(H^{1}(\Tilde{C}, \cO_{\Tilde{C}}))} = \coker(F_{\Tilde{C}}).$$ 
For each component $C_i$ of $\Tilde{C}$, there is a Hasse-Witt triple $(H^1(C_i, \cO_{C_i}), F_{C_i}, \Psi_{C_i})$, and $\Psi_{C_i}$ is unique up to isomorphism by Theorem \ref{thm:moonen}. The collection $(\Psi_{C_i})_i$ defines 
$\Psi_{\Tilde{C}}$ uniquely up to isomorphism. Consequently, the~Dieudonn\'e module associated with the Hasse-Witt triple $(H^1(C, \cO_C), F_C, \Psi_C)$ is isomorphic to $\mathbb{D}_{l, \text{ord}}\oplus \bigoplus_{i = 1}^n \mathbb{D}_{C_i}$, where $\mathbb{D}_{C_i} \cong \mathbb{D}_{\cJ_{C_i}}$ is the Dieudonn\'e module of $C_i$. Furthermore, it holds that $\mu(C) = \mu(\cJ_{C}) = \mu(\Tilde{C})$. 
\end{proof}

\begin{rem}
Let $C$ be a stable curve with $m$ singular points and $C_i$, $1\leq i \leq n$ its irreducible components. In general, we are not aware of a straightforward combinatorial description of $\mu(C)$ in terms of $\mu(C_i)$, for $1 \leq i \leq n$. However, the collection of $(\mu(C_i))_{1\leq i\leq n}$ (together with $m$) uniquely determines $\mu(C)$ and one can explicitly compute it using \cite[~9.1]{oort}.
\label{rem:muCnorm}
\end{rem}

\section{Inductive results}
\label{sec:inductiveres}

Throughout this section, we use that $\cMbar_g$ and $\cA_g$ are smooth stacks, which implies that the codimension of the intersection of any two of their closed subvarieties is at most the sum of the codimensions. Recall $\dim \cMbar_g = 3g - 3$ for $g\geq 2$ and $\dim \cA_g = \frac{g(g + 1)}{2}$ for $g \geq 1$. Let $\cM_g^{ct} = \cMbar_g - \Delta_0$ denote the moduli space of genus-$g$ curves of compact type, and note that $j(\cM_g^{ct})\subseteq \cA_g$.
For $g$ and $n$ such that $2g - 2 + n > 0$, we denote by $\cMbar_{g, n}$ the moduli space of $n$-pointed stable curves of genus $g$. 
The clutching morphisms $\kappa_{1, 1}$ and $\kappa_2$ as defined in \cite[Section~3]{knudsen} are finite and given by 
$$\kappa_{1, 1}: \cMbar_{g, 1} \times \cMbar_{g', 1} \to \cMbar_{g + g'}\quad \text{and} \quad \kappa_{2}: \cMbar_{g - 1, 2}  \to \cMbar_{g}, $$ 
where $\kappa_{1, 1}$ identifies the labeled points on stable genus-$g$ and genus-$g'$ curves, and $\kappa_2$ identifies the labeled points on a stable genus-$(g - 1)$ curve.
For any $g\geq 2$ and $n\geq 1$, let~$Z_{\mu}(\cMbar_{g, n})$ denote the inverse image of $Z_{\mu}(\cMbar_g)$ under the forgetful morphism ${\cMbar_{g, n} \to \cMbar_g}$ (whose fibers are $n$-dimensional). For $(g, n) = (1, 1)$, the locus $Z_{\o}(\cMbar_{1, 1})$ represents the  
locus of ordinary stable elliptic curves, while $Z_{[1]}(\cMbar_{1, 1})$ represents the locus of supersingular elliptic curves.  Lastly, let $Z_{\mu}(\cM_{g}) = Z_{\mu}(\cMbar_g) \cap \cM_g$.

 The dimensions and many other geometric properties of the Ekedahl-Oort strata of $\cA_g$ are known, including those of the $\prank \leq f$ and $a$-number $\geq a$ strata; see Section \ref{sec:eo_types} (or \cite{oort} and \cite{vdgeercycle}). However, our understanding remains limited for $\cMbar_g$, $g\geq~4$.
Let $$V_f(\cMbar_{g}) = \bZ_{[g - f]}(\cMbar_{g}) \text{ and } T_a(\cMbar_{g}) = \bZ_{[a, a-1, \ldots, 2, 1]}(\cMbar_{g})$$ denote the loci of stable curves of genus $g$ with $\prank$ at most~$f$ and $a$-number at least~$a$, respectively. Note that $\bZ_{[g - f, 1]}(\cMbar_{g})$ equals the sublocus of $V_f(\cMbar_g)$ with $a$-number at least $2$, if $f \leq g - 2$. The main results for $\cMbar_g$ with $g \geq 2$ are summarized in the following theorem; additionally, see Remark \ref{rem:priessmootheo}.

\begin{thm} 
\label{thm:fvgdp}
The following results hold. 
\begin{enumerate}
    \item \cite[Theorem 2.3]{fabervdgeer} For any $g\geq 2$ and $0\leq f \leq g$, the locus $V_f(\cMbar_{g})$ is pure of codimension $g - f$ in $\cMbar_g$.
    \item \cite[Proposition 3.7]{pries_a_number} Let $g_0 \geq 2$ and $0\leq f_0 \leq g_0$. If every generic point of $V_{f_0}(\cMbar_{g_0})$ has $a$-number $1$ then for any $g> g_0$ and $f = g - (g_0 - f_0)>0$, any generic point of $V_{f}(\cMbar_{g})$ has $a$-number $1$.
\end{enumerate}
\end{thm}

The main technique used in proving the previous results is an inductive argument based on our understanding of $\cMbar_{g - 1}$. This approach focuses on loci whose intersection with the boundary divisor $\Delta_0 = \kappa_2(\cMbar_{g - 1, 2})$ is non-empty, as described in the following remark.

\begin{rem}[Diaz-Looijenga bound]
   If $\Gamma$ is a closed subvariety of $\cMbar_g$ whose codimension in $\cMbar_g$ is at most $g-1$, then $\Gamma \cap \Delta_0 \neq 0$. This was first observed (in characteristic~$0$) by Diaz in \cite[page~80]{diaz}  and follows from the Diaz-Looijenga bound in \cite[page~412]{looijenga}; see also \cite[Lemma~2.4]{fabervdgeer}.
    \label{rem:diazbound}
\end{rem}

Theorem \ref{thm:fvgdp} $(1)$ can be understood informally as a consequence of the fact that $V_0(\cMbar_g)$ has codimension $g$ in $\cMbar_g$ for all $g\geq 2$. Moreover, the knowledge about $V_0(\cMbar_2)$ and $V_0(\cMbar_3)$ (see Remark \ref{rem:priessmootheo} below) together with Theorem \ref{thm:fvgdp}~$(2)$ leads to the result below. This result is essentially derived in the proofs of \cite[Corollary~4.5]{pries_a_number} and \cite[Proposition~4.9]{pries_a_number}.

\begin{cor}[{\cite{pries_a_number}}] 
\label{cor:pries_anum}
For any $g\geq 3$ and $p > 0$, the following results hold. 
\begin{enumerate}
   \item The locus of curves in $V_{g - 2}(\cMbar_g)$ with $a$-number $\geq 2$ is pure of codimension $1$ in $V_{g - 2}(\cMbar_g)$, and it contains a point that represents a smooth curve.
    \item The locus of curves in $V_{g - 3}(\cMbar_g)$ with $a$-number $\geq 2$ is pure of codimension $1$ in $V_{g - 3}(\cMbar_g)$, and it contains a point that represents a smooth curve.
\end{enumerate}
\end{cor}

Our goal in this section is to show the full potential of the inductive technique introduced in \cite{fabervdgeer}, used to prove both statements mentioned in Theorem \ref{thm:fvgdp}.

\subsection{Inductive results about the Ekedahl-Oort strata}

Fix $g \geq 2$ and let $\mu$ and $\mu'$ be two Young diagrams. The inclusion $Z_{\mu'} \subseteq \bZ_{\mu}$ of the Ekedahl-Oort loci in $\cA_g$ can be determined by the combinatorial description in \eqref{eqn:shuffle}. This description involves specific relations between elements in the Weyl group $W_g$ and shuffles. According to \cite[page 599]{ekedahlvdgeer}, there exists an injective group homomorphism $\rho: W_g \to W_{g+1}$ that preserves both the Bruhat-Chevalley order and the shuffling. Consequently, this homomorphism respects \eqref{eqn:shuffle}, ensuring that if $Z_{\mu'} \subseteq \bZ_{\mu}$ in $\cA_g$, then $Z_{\mu'} \subseteq \bZ_{\mu}$ also holds in $\cA_{g+1}$.

Let $\mu' = [\mu_1', \ldots, \mu_{n'}']$ and consider the following condition on the Young diagram $\mu$: 
\begin{equation}
    Z_{\mu'} \subseteq \bZ_{\mu} \text{ in } \cA_{g + 1} \implies Z_{\mu'} \subseteq \bZ_{\mu} \text{ in } \cA_{g} \text{ for each }\mu' \text{ with } \mu_1' < g + 1.
\label{eqn:technical_cond}
\end{equation}
Note that \eqref{eqn:shuffle} depends on $g$ and the structure of $W_g$. The condition $Z_{\mu'} \subseteq \bZ_{\mu}$ in~$\cA_{g + 1}$ is governed by a relation in $W_{g + 1}$ involving shuffles, while in $\cA_g$ it relates to elements of~$W_g$. As a result, \eqref{eqn:technical_cond} is not trivial, and to present Theorem \ref{thm:main1} in full generality, we must focus on cases where this condition is met. Below, we mention some instances.

\begin{rem}[Young diagrams describing $\prank \leq f$ and $a$-number $\geq a$ loci satisfy~\eqref{eqn:technical_cond}]
\label{rem:ordering}
Note that   \eqref{eqn:technical_cond} is satisfied for all $\mu$ such that $\bZ_{\mu} = \cup_{\mu'\leq \mu} Z_{\mu'}$ in $\cA_{g + 1}$ with $g\geq \mu_1$; see \eqref{eq:finaltypes_EOclosure}.
Using \eqref{eqn:oort125}, for any $a, f \geq 0$ and $g \geq 2$ such that $1\leq a + f \leq g$, we see that the Young diagram 
\begin{equation}
\label{eqn:mu_prank_anum_great}
\mu = [g - f, \underbrace{a - 1, a - 2, \ldots, 2, 1}_{a-1}]   
\end{equation}
satisfies $\bZ_{\mu} = \cup_{\mu'\leq \mu} Z_{\mu'}$ in $\cA_{g}$.  
This follows from $\bZ_{\mu} = V_f\cap T_a$ being the closed locus of principally polarized abelian varieties with $\prank \leq f$ and $a$-number $\geq a$; see \cite[Corollary 1.5]{oort_subv} and \cite[page~44]{koblitz}. Thus, for $\mu$ as in $\eqref{eqn:mu_prank_anum_great}$, the technical condition \eqref{eqn:technical_cond} is automatically met for any $g$. We will use this observation below without further comment. 
\end{rem}

\begin{exmp} 
\label{exmp:eo}
Let us present a few examples of Young diagrams that satisfy \eqref{eqn:technical_cond} or \eqref{eqn:mu_prank_anum_great}.
\begin{itemize}
    \item As noted earlier, the types $\mu = [g - f]$ and $\mu = [g - f, 1]$ define the locus of stable genus-$g$ with $\prank\leq f$,  and its sublocus where the $a$-number is at least $2$, respectively. Both types are of the form \eqref{eqn:mu_prank_anum_great} and thus satisfy \eqref{eqn:technical_cond} for any~$g$.
    \item Similarly, $\mu = [3, 2, 1]$ is of the form \eqref{eqn:mu_prank_anum_great} and therefore satisfies \eqref{eqn:technical_cond} for any $g$. Recall that $\bZ_{[3, 2, 1]}(\cMbar_{g})$ is the locus of stable curves of genus $g$ with $a$-number $\geq 3$.
    \item The Young diagram $\mu = [3, 2]$ satisfies \eqref{eqn:technical_cond} at least for $g = 3, 4$, even though it is not of the form \eqref{eqn:mu_prank_anum_great}. This is immediate for $g = 3$, and follows from Example~\ref{exmp:4132} for $g = 4$. 
\end{itemize}
\end{exmp}

Having addressed the technical conditions, we now turn to our second main result. It relies on the intrinsic nature of Definition \ref{dfn:eo_stable} and Theorem \ref{thm:muCnorm}, and shows the full potential of the technique used in the proof of \cite[Theorem 2.3]{fabervdgeer}. 

\begin{thm}
Let $g\geq 2$, let $\mu$ be a Young diagram that satisfies the condition~\eqref{eqn:technical_cond}, and let $d = 3g - 3 - \max \{\dim\Gamma':\Gamma' \text{ a component of }\bZ_{\mu}(\cMbar_g)\}.$ If $\Gamma\subseteq \bZ_{\mu}(\cMbar_{g + 1})$ is a component such that $\Gamma \cap \Delta_0 \neq \emptyset$, then the codimension of $\Gamma$ in $\cMbar_{g + 1}$ is at least~$d$.
\label{thm:main1}
\end{thm}
{\begin{proof}
Let $C$ be any stable curve that corresponds to some point of $\Gamma\cap \Delta_0 \neq \emptyset$. Let $\mu' = [\mu_1', \ldots, \mu_{n'}']$ denote its Ekedahl-Oort type, and observe that $Z_{\mu'} \subseteq \bZ_{\mu}$ in~$\cA_{g + 1}$ by \eqref{eqn:oort125}. Since $f(C) = g + 1 - \mu_1'$ and the  $\prank$ of any curve lying in $\Delta_0 = \kappa_2(\cMbar_{g, 2})$ has to be at least~$1$, it follows that  $\mu'_1 < g + 1$.

Drawing on the insight from Theorem \ref{thm:muCnorm}, the normalization of any singular curve $C$ with Ekedahl-Oort type $\mu'$ also has Ekedahl-Oort type $\mu'$. Since $\mu$ satisfies the assumption \eqref{eqn:technical_cond} and we have that $Z_{\mu'} \subseteq \bZ_{\mu} \text{ in } \cA_{g + 1}$ as shown above, we conclude that $Z_{\mu'}\subseteq\bZ_{\mu}$ in $\cA_g$. As a result, 
\begin{equation}
\Gamma\cap \Delta_0\subseteq \kappa_2(\bZ_{\mu}(\cMbar_{g, 2})).
\label{eqn:intersection_contained_in_Delta0}
\end{equation}
Note that $\dim \kappa_2(\bZ_{\mu}(\cMbar_{g, 2})) = \dim \bZ_{\mu}(\cMbar_{g}) + 2$. The boundary component $\Delta_0$ has codimension $1$ in $\cMbar_g$. Since $\cMbar_g$ is a smooth stack, it follows that every component of $\Gamma \cap \Delta_0$ has codimension at most $1$ in $\Gamma$. Therefore, \eqref{eqn:intersection_contained_in_Delta0} implies that  $$\dim \Gamma  \leq \dim \kappa_2(\bZ_{\mu}(\cMbar_{g, 2})) + 1 \leq (3g - 3 -d + 2) + 1 = 3(g + 1) - 3  - d.$$ 
This establishes the theorem since saying that the codimension of $\Gamma$ in $\cMbar_{g + 1}$ is at least $d$ is the same as showing that $\dim \Gamma \leq 3(g + 1) - 3 - d$.
\end{proof}}

Let us first comment on how Theorem \ref{thm:main1} gives a proof of Theorem \ref{thm:fvgdp}~$(1)$ and $(2)$. 
For the part $(1)$, the upper bound on $\dim V_f (\cMbar_g)$ follows from Theorem~\ref{thm:main1} for the choice $\mu = [g - f]$ and induction on $g$ with $g - f$ fixed since $ V_f (\cMbar_g) = \bZ_{[g - f]}(\cMbar_{g})$, while the lower bound  follows from a purity result (\cite[Lemma~1.6]{oort_subv}, see also \cite[Lemma~2.1]{fabervdgeer}). Theorem~\ref{thm:fvgdp} $(2)$ is equivalent to showing that the sublocus $\bZ_{[g -f , 1]}(\cMbar_{g})$ of $V_f(\cMbar_{g})$ consisting of the curves whose $a$-number is at least $2$ satisfies $\dim \bZ_{[g -f , 1]}(\cMbar_{g}) < \dim V_f(\cMbar_{g})$; see Example~\ref{exmp:eo}. This follows from Theorem \ref{thm:main1} for the choice $\mu = [g - f, 1]$ and induction on~$g$ with~$g - f$ fixed. Corollary~\ref{cor:pries_anum} follows by induction on $g$ since we know that $\bZ_{[2, 1]}(\cM_2)$ (resp. $\bZ_{[3, 1]}(\cM_3)$) is non-empty and its closure has codimension $3$ in~$\cMbar_2$ (resp. $4$ in $\cMbar_3$). 

\begin{rem} 
\label{rem:codim_d_general}
Here is some additional commentary on the preceding theorem.
\begin{enumerate}
    \item The condition $\Gamma \cap \Delta_0 \neq \emptyset$ in Theorem \ref{thm:main1} is \textit{automatically} satisfied if $g \geq \sum \mu_i$ by the Diaz-Looijenga bound; see Remark \ref{rem:diazbound}.
    \item A chain of $g$ supersingular elliptic curves is a superspecial stable curve, i.e., a stable curve with Ekedahl-Oort type $[g, g-1, \ldots, 2, 1]$. Hence, the loci $$Z_{[g, g-1, \ldots, 2, 1]}(\cMbar_{g}) \subseteq \bZ_{\mu}(\cMbar_g)$$ are always non-empty, though $Z_{\mu}(\cMbar_g)$ or $\bZ_{\mu}(\cM_{g})$ might be empty. 
    \item For the choice of~$\mu$ as in \eqref{eqn:mu_prank_anum_great}, the preceding theorem provides an upper bound on the dimension of loci in $\cMbar_{g+1}$ where the $p$-rank is $\leq f+1$ and the $a$-number is $\geq a$, in terms of the bound on the dimension of loci in $\cMbar_g$ where the $p$-rank is $\leq f$ and the $a$-number is $\geq a$.
\end{enumerate}
\end{rem}

\begin{rem}
In \cite[Theorem 6.4]{pries_current_results}, Pries  shows that if a component $\Gamma_0\subseteq Z_{\mu}(\cM_{g_0})$ of codimension $d = \sum \mu_i$ in $\cM_{g_0}$ exists, then there is a component $\Gamma$ of $Z_{\mu}(\cM_{g})$ of codimension $d = \sum \mu_i$ in $\cM_g$ for any $g\geq g_0$. (Both $\Gamma_0$ and $\Gamma$ consist of isomorphism classes of smooth curves.) By \cite[3.3.1, 3.3.2, and Theorem 6.4]{pries_current_results}, given a $\mu$, there is a non-empty component of $Z_{\mu}(\cM_{g})$ of the expected dimension in characteristic $p$ in the following cases:
\begin{enumerate}
    \item $\mu  \in \{[1], [2], [3], [3, 1]\}$ for any $p>0$, 
    \item $\mu = [2, 1]$, for $g = 2$, $p>3$  or $g \geq 3$, $p>0$, or
    \item $\mu  \in \{[3, 2], [3, 2, 1]\}$ and $p>2$.
\end{enumerate} 
\label{rem:priessmootheo}
See also \cite[Theorem 5.12]{oort_hess}, especially for the case $g = 3$ and $\mu = [3, 2, 1]$.

Consider the differences between \cite[Theorem 6.4]{pries_current_results} and Theorem \ref{thm:main1} and how they complement each other under the assumption $g > \sum \mu_i.$ 
By the construction in \cite[Theorem 6.4]{pries_current_results}, each $\Gamma_0$ as above yields a $\Gamma$ in $\cMbar_g$ whose closure 
contains the image of~$\cMbar_{1, 1}\times \Gamma_0'$ under the clutching morphism $\kappa_{1, 1}$, where~$\Gamma_0'$ is the pullback of~$\Gamma_0$ under the forgetful morphism $\cMbar_{g_0, 1} \to \cMbar_{g_0}$. This does not necessarily describe all components of $Z_{\mu}(\cM_{g})$. Conversely, while Theorem \ref{thm:main1} does not guarantee the existence of smooth curves with prescribed Ekedahl-Oort type $\mu$, it provides an \textit{upper} bound of \textit{all} components $\Gamma$ of $\bZ_{\mu}(\cM_{g})$; recall that $\overline{\Gamma} \cap \Delta_0 \neq \emptyset$ whenever $g > \sum \mu_i$. Once we know there is one such curve (e.g., with the help of \cite[Theorem~6.4]{pries_current_results}), we know the \textit{lower} bound $\dim \Gamma\geq 3g - 3 - \sum \mu_i$ for all components $\Gamma$ of $\bZ_{\mu}(\cM_{g})$ using that $\cA_g$ is a smooth stack and $\dim \Gamma = \dim j(\Gamma)$. These combined findings can provide the exact dimension for any such $\Gamma$.
\end{rem}

\section{Some applications}
\label{sec:app}

Let us illustrate how we can obtain new results about the dimensions of Ekedahl-Oort strata $\bZ_{\mu}(\cMbar_g)$, by looking at $\mu = [3, 2]$ and $\mu = [3, 2, 1]$. As noted in Example~\ref{exmp:eo}, the type $\mu = [3, 2]$ does not define a $\prank \leq f$ and $a$-number $\geq a$ locus, but it satisfies~\eqref{eqn:technical_cond} for $g = 3$ and $g = 4$.
 In the case $\mu = [3, 2]$,  we see the first example where the dimensions of components in $\bZ_{\mu}(\cMbar_g)$ for $g \geq 3$ do not match the expected $3g - 3 - \sum \mu_i$. Namely,~$\bZ_{[3, 2]}(\cMbar_{3})$ contains $\kappa_{1, 1}(\bZ_{[1]}(\cMbar_{1, 1}) \times \bZ_{[2]}(\cMbar_{2, 1}))$, which is a $2$-dimensional locus. 

\begin{prop}
\label{prop:32_genus45}
The following results hold. 
\begin{itemize}
    \item Any component $\Gamma$ of $\bZ_{[3, 2]}(\cMbar_{4})$ satisfies $\dim \Gamma \leq 5$.
    \item Any component $\Gamma$ of $\bZ_{[3, 2]}(\cMbar_{5})$ satisfies $\dim \Gamma \leq 8$.
\end{itemize}
\end{prop}

\begin{proof}
 First, note that ${\dim \bZ_{[3, 2]}(\cMbar_{3})\leq 2}$. Namely, $\cM_3\cap  \bZ_{[3, 2]}(\cMbar_{3})$ has dimension $\leq 1$ by \cite[Corollary~11.2]{oort}, while  $(\cMbar_3 - \cM_3) \cap \bZ_{[3, 2]}(\cMbar_{3})$ equals $\kappa_{1, 1}(\bZ_{[1]}(\cMbar_{1, 1}) \times \bZ_{[2]}(\cMbar_{2, 1}))$ and thus has dimension $2$, which results in $\dim \bZ_{[3, 2]}(\cMbar_{3})\leq 2$. 
 
 Let $\Gamma$ be a component of $\bZ_{[3, 2]}(\cMbar_4)$ in characteristic $p>0$. Hence, if ${\Gamma}\cap \Delta_0 \neq \emptyset$, we compute $\dim \Gamma \leq \dim \kappa_{2}(\bZ_{[3, 2]}(\cMbar_{3, 2})) + 1 \leq 5$
using Theorem \ref{thm:main1} and $\dim \bZ_{[3, 2]}(\cMbar_{3})\leq 2$.
Otherwise, $\Gamma \subseteq \cM^{ct}_4$, and $\dim \Gamma \leq 5$ by the the Diaz-Looijenga bound.

Now, let $\Gamma$ be a component of $\bZ_{[3, 2]}(\cMbar_{5})$. Similarly as above, using Theorem~\ref{thm:main1} and $\dim \bZ_{[3, 2]}(\cMbar_{4})\leq 5$, we conclude $\dim \Gamma \leq 8$ if ${\Gamma}\cap \Delta_0 \neq \emptyset$. Otherwise,  $\dim \Gamma \leq 7$ by the Diaz-Looijenga bound.
\end{proof}

We will use the following lemma in some of the descriptions below. 

\begin{lem}
\label{lem:intersect_delta0}
In any characteristic $p>0$, every component of the Ekedahl-Oort loci $\bZ_{[1]}(\cMbar_{2})$, $\bZ_{[2]}(\cMbar_{3})$, and $\bZ_{[2, 1]}(\cMbar_{3})$ has a non-empty intersection with $\Delta_0$.   
\end{lem}
\begin{proof}
First, $\Gamma = \bZ_{[1]}(\cM_2)$ is irreducible since $\overline{j(\Gamma)} = \bZ_{[1]}$ is irreducible in $\cA_2$ by \cite[Theorem 11.5]{ekedahlvdgeer}; in fact,  $\overline{j(\Gamma)} = V_1 \subseteq \cA_2$ so its irreducibility was observed even before. Furthermore, we have $\overline{\Gamma} = \bZ_{[1]}(\cMbar_{2})$ by the dimension count; $\bZ_{[1]}(\cMbar_{2}) \cap \Delta_0 \neq \emptyset$ follows from the Diaz-Looijenga bound in Remark \ref{rem:diazbound}. A similar argument gives a proof for $\bZ_{[2]}(\cMbar_{3})$; see also \cite[Lemma 3.3]{achterpries}.

In the case of $\bZ_{[2, 1]}(\cMbar_{3})$, we cannot use the Diaz-Looijenga bound to get the result. Consider the locus $\bZ_{[2, 1]}\subseteq\Tilde{\cA}_3.$ By \cite[Theorem 11.5]{ekedahlvdgeer}, this locus is irreducible. Moreover, by the dimension count, its generic point represents a Jacobian of a smooth curve. These two facts imply that $\bZ_{[2,1]}(\cMbar_3)$ has a unique component $\Gamma$, whose generic point corresponds to a smooth curve, and that ${j(\Gamma)} = \bZ_{[2,1]} \subseteq \Tilde{\cA}_3.$ Since $\bZ_{[2,1]}$ contains a point that represents a generalized Jacobian of a curve in~$\Delta_0 \subseteq \cMbar_3$, it follows that $\Gamma \cap \Delta_0 \neq \emptyset$. 
The remaining $3$-dimensional components of $\bZ_{[2,1]}(\cMbar_3)$ that are not entirely contained in~$\Delta_0$ lie in 
$$\kappa_{1,1}(\bZ_{[1]}(\cMbar_{2,1})\times \bZ_{[1]}(\cMbar_{1,1})).$$ 
By the first part of the proof, they have a non-empty intersection with $\Delta_0$. This finishes the proof.
\end{proof}

Now, we consider the loci $\bZ_{[3, 2, 1]}(\cMbar_{g})$ for any $g\geq 3$; see Example \ref{exmp:eo}. 

\begin{prop}
\label{prop:321}
Let $\Gamma$ be a component of $ \bZ_{[3, 2, 1]}(\cMbar_{g})$, for some $g \geq 3$. Then, the dimension of $\Gamma$ is at most $3g - 8$.    
\end{prop}

\begin{proof}
First, we consider the locus $\bZ_{[3, 2, 1]}(\cMbar_{3})$ of stable superspecial curves of genus $3$. Note that it contains the $1$-dimensional locus $\kappa_{1, 1}(\bZ_{[1]}(\cMbar_{1, 1})\times \bZ_{[2, 1]}(\cMbar_{2, 1}))$, while its sublocus representing smooth curves is finite if non-empty. Therefore, $\dim \bZ_{[3, 2, 1]}(\cMbar_{3}) \leq 1.$

Let $\Gamma$ be a component of $\bZ_{[3, 2, 1]}(\cMbar_{4})$. If $\Gamma \cap \Delta_0 \neq \emptyset$, Theorem \ref{thm:main1} and the observation we just made imply $$\dim \Gamma \leq \dim \kappa_{2}(\bZ_{[3, 2, 1]}(\cMbar_{3, 2})) + 1 \leq 4.$$ Otherwise, $\Gamma \subseteq \cM^{ct}_4$. If the generic point of $\Gamma$ corresponds to a smooth curve, we conclude $\dim \Gamma = \dim j(\Gamma) \leq 4 = \dim\bZ_{[3, 2, 1]}$. If not, $\Gamma \subseteq \Delta_1 \cup \Delta_2$ and $\dim \Gamma \leq 3$. Namely, an irreducible component $\Gamma$ of $\bZ{[3, 2, 1]}(\cMbar_4) \cap \Delta$ with $\dim \Gamma \geq 3$ is either a component $\Gamma_1 \subseteq \kappa_{1, 1}(\bZ_{[1]}(\cMbar_{1, 1}) \times \bZ_{[2, 1]}(\cMbar_{3, 1}))$, a component $\Gamma_2 \subseteq \kappa_{1, 1}(\bZ_{[1]}(\cMbar_{2, 1}) \times \bZ_{[2, 1]}(\cMbar_{2, 1}))$, or a component $\Gamma_3 \subseteq \kappa_{1, 1}(\bZ_{[2]}(\cMbar_{2, 1}) \times \bZ_{[2, 1]}(\cMbar_{2, 1}))$. Note that $\dim \Gamma_1 = \dim \Gamma_3 = 3$ and $\dim \Gamma_2 = 4$, while Lemma \ref{lem:intersect_delta0} implies that $\Gamma_1 \cap \Delta_0 \neq \emptyset$ and $\Gamma_2 \cap \Delta_0 \neq \emptyset$.

Next, we obtain $$\dim \bZ_{[3, 2, 1]}(\cMbar_{5}) \leq 7 \quad\text{ and}\quad \dim \bZ_{[3, 2, 1]}(\cMbar_{6}) \leq 10.$$ For any component~$\Gamma$ of either $\bZ_{[3, 2, 1]}(\cMbar_{5})$ or $\bZ_{[3, 2, 1]}(\cMbar_{6})$, the claimed upper bound on $\dim \Gamma$ follows from the Diaz-Looijenga bound if $\Gamma \cap \Delta_0 = \emptyset$ or from Theorem \ref{thm:main1} if $\Gamma \cap \Delta_0 \neq \emptyset$.

Finally, let $g \geq 7$ and let $\Gamma$ be a component of $\bZ_{[3, 2, 1]}(\cMbar_{g})$. By the Diaz-Looijenga bound, $\Gamma \cap \Delta_0  = \emptyset$ implies that $\dim \Gamma \leq 3g - 3 - g < 3g - 8$. Otherwise, $\Gamma \cap \Delta_0 \neq 0$ so $\dim \Gamma \leq 3g - 8$ by Theorem \ref{thm:main1} and induction on $g$.
\end{proof} 

\begin{rem}
Although the 1-dimensional locus $\kappa_{1, 1}(\bZ_{[1]}(\cMbar_{1, 1}) \times \bZ_{[2, 1]}(\cMbar_{2, 1}))$ consists entirely of singular curves, its existence suggests the possible existence of a 4-dimensional family of smooth genus-$4$ curves with Ekedahl-Oort type $[3, 2, 1]$, which contrasts with the expected dimension of $3 = \dim \cM_4 - (3 + 2 + 1)$ for the components of $Z_{[3, 2, 1]}(\cM_4)$. However, for example, this is not the case for~$p = 2$, as $Z_{[3,2,1]}(\cM_4) = \emptyset$ in characteristic $2$ (see Proposition \ref{prop:eo_conclusion_char2} below). 
\end{rem}

Finally, let $T_a(\cMbar_g) = \bZ_{[a, a-1, \ldots, 2, 1]}(\cMbar_g)$ denote the locus of stable curves of genus $g$ with $a$-number $\geq a$. We summarize some of our findings below.  
\begin{cor} Let $g \geq 3$ and let $p>0$ be any prime number. Then, the following results hold in characteristic $p$: $\dim T_{1}(\cMbar_g)  = 3g - 4$, $\dim T_{2}(\cMbar_g)  = 3g - 6$, and $\dim T_{3}(\cMbar_g)  \leq 3g - 8$.
\end{cor}
\begin{proof}
The conclusions for $T_{1}(\cMbar_g)  = V_{g - 1}(\cMbar_g)$ and  $T_2(\cMbar_g)$ follow from Theorem~\ref{thm:fvgdp} (1) and Corollary \ref{cor:pries_anum} (1), respectively. Lastly, the conclusion  $$\dim T_{3}(\cMbar_g)  \leq 3g - 8$$ follows from Proposition \ref{prop:321} since $T_{3}(\cMbar_g) = \bZ_{[3, 2, 1]}(\cMbar_g)$.  
\end{proof}

\subsection{On Ekedahl-Oort loci of smooth hyperelliptic curves of genus $4$}

Consider the (closed) locus  $\overline{\cH}_4$ of stable hyperelliptic curves of genus~$4$, which is of codimension $4 - 2 = 2$ in $\cMbar_4$. Let  $Z_{[4, 3]}(\cH_{4})$ denote the locus of smooth hyperelliptic curves $C$ of genus $4$ over $k$ with Ekedahl-Oort type~$[4, 3]$, or equivalently, with
\begin{equation}
\rank_k(F_C) = 2 \quad \text{ and }\quad \rank_k(F_C^2) = 0.
\label{eqn:43_eqcond}
\end{equation}

We have the following theorem.

\begin{thm}
\label{thm:he_implies_eo}
Let $\mu$ be a Young diagram such that ${\mu > [4, 3]}$. 
If $Z_{[4, 3]}(\cH_{4})$ is non-empty and finite in characteristic $p>0$, then $Z_{\mu}(\cM_4)$ is pure of codimension $cd(\mu)$ in $\cM_4$.
\end{thm}

\begin{proof}
Let $\mu $ be a Young diagram with $\mu > [4, 3]$. By \cite[Theorem 11.5]{ekedahlvdgeer}, $\bZ_{\mu}\subseteq \cA_4$ is irreducible. 
For any component $\Gamma^0\subseteq Z_{\mu}(\cM_4)$ we have $$\dim j(\Gamma^0) = \dim \Gamma^0 = \dim \overline{\Gamma^0}.$$ 
Since $\cJ_4 = j(\cM_4^{ct})$ has codimension $1$ in $\cA_4$ and $\cA_4$ is smooth, 
$\bZ_{\mu}\cap \cJ_4$ is either pure of dimension $9 - \sum_{i = 1}^n\mu_i$, or it is pure of dimension $10 - \sum_{i = 1}^n\mu_i$; in the latter case $\bZ_{\mu} \subseteq \cJ_4$.

Assume the latter is true and let $\Gamma = \overline{\Gamma^0}$ be a component of $\bZ_{\mu}(\cMbar_4)$ containing a point~$[C]$ of $Z_{[4, 3]}(\cH_{4})$ that represents a smooth hyperelliptic curve $C$ with Ekedahl-Oort type $\mu(C) = [4, 3]$; note that $j(\Gamma\cap \cM_g^{ct}) = \bZ_{\mu} \subseteq \cA_4$. Since $\bZ_{[4, 3]} \subseteq \bZ_{\mu}$ and $\dim \bZ_{[4, 3]} = 3$, there is a closed irreducible locus $\Gamma' \subseteq \bZ_{[4, 3]}(\cMbar_{4})$ such that $$\{[C]\} \subseteq \Gamma' \subseteq \Gamma$$ and $\dim \Gamma' = 3 = \dim j(\Gamma')$. Now, $\Gamma' \cap \overline{\cH}_{4} \neq \emptyset$ and using that $\cMbar_4$ is a smooth stack, we find that every component of $\Gamma' \cap \overline{\cH}_{4} \neq \emptyset$ has dimension at least  $3 - 2 = 1$. However, $\{[C]\}$ is a $0$-dimensional component of $\Gamma' \cap \overline{\cH}_{4}$. This contradicts our assumption and we conclude that $\dim Z_{\mu}(\cM_4) = 9  - cd(\mu)$ if $Z_{\mu}(\cM_4) \neq \emptyset$.

To see that $Z_{\mu}(\cM_4) \neq \o$, we consider the component $\Gamma \subseteq \bZ_{\mu}(\cMbar_4)$ such that $\{[C]\} \subseteq \Gamma$ as above, and observe that the dimension count gives us that the generic point of $\Gamma$ has to correspond to a smooth curve with Ekedahl-Oort type $\mu$.
\end{proof}

\begin{rem}
\label{rem:he43finite_implies_ss3dim}
A variation of the preceding criterion can be applied to compute dimensions of certain Newton polygon strata. Namely, assuming that $Z_{[4, 3]}(\cH_{4})$ is non-empty and finite, and noting that $Z_{[4, 3]} \subseteq \cS_4$, where $\cS_4$ is the supersingular locus of $\cA_4$, the same argument results in a component $W$ of the supersingular locus of $\cM_4$ with $\dim W = 3$ in characteristic~$p$. Furthermore, it follows that the codimension of every component of a non-supersingular Newton polygon locus in $\cM_4$ equals the expected value. See also \cite[Remark 5.3]{dd2}.
\end{rem}

\section{Examples in characteristic $p = 2$ and $p = 3$}
\label{sec:exm}

In this section, we explore the Ekedahl-Oort stratification of moduli spaces of smooth genus-$g$ curves for $g \geq 4$ in characteristics $p = 2$ and $p = 3$.

\subsection{Ekedahl-Oort stratification of $\cM_4$ in characteristic $2$}
 
By Remark \ref{rem:priessmootheo}, there is a smooth curve of genus $4$ with Ekedahl-Oort type $\mu < [3, 2]$ in any characteristic $p > 0$, while there exists a smooth curve of genus $4$ with Ekedahl-Oort type $\mu \leq [3, 2]$ in characteristic $p > 2$.  We consider the remaining case $\mu = [3, 2]$ and $p = 2$ below. 

\begin{lem} 
\label{lem:char2_eo32_genus4}
In characteristic $p = 2$, there is a smooth genus $4$ curve $C$ with $\mu(C) = [3, 2]$.   
\end{lem}
\begin{proof}
From Xarles' list of all genus-$4$ curves defined over $\F_2$ (up to $\F_2$-isomorphism) in \cite{xarles}, we identify the hyperelliptic curve $$C: y^2 + x^2y = x^9 + x$$ as a potential candidate. After the change of coordinates, $C: y^2 - y = \sum_{\alpha \in \{0, 1\} } f_{\alpha}(x_{\alpha})$ is in the form of  \cite[Notation~1.1]{elkinpries}, for $\alpha \in \{0, 1\}$ with $x_{\alpha}  = (x - \alpha)^{-1}$ and $f_{\alpha}(x) \in \F_2[x]$ of degree $\deg f_{\alpha} = d_{\alpha}$, where $d_{0} = 3$ and $d_{1} = 5$. We compute $\mu(C) = [3, 2]$ using \cite[Theorem 5.2]{elkinpries}.
\end{proof}

A similar proof to that of \cite[Corollary 6.6]{dd2} establishes a slight generalization below.

\begin{prop}
\label{prop:eo_conclusion_char2}
In characteristic $p = 2$, for any Young diagram $\mu = [\mu_1, \ldots, \mu_n]$ with $\mu > [4, 3]$, $Z_{\mu}(\cM_4)$ is non-empty and is pure of codimension $\sum \mu_i$ in $\cM_4$. For any other Young diagram $\mu$, namely $\mu = [4, 3]$ or $\mu \leq [3, 2, 1]$, $Z_{\mu}(\cM_4)$ is empty.
\end{prop}

\begin{proof}
In characteristic $2$, the (non-)emptiness of $Z_{\mu}(\cM_4)$ given a Young diagrams $\mu$ as above with $\mu_1 = 4$ is as stated by \cite[Corollary 6.6]{dd2}, for $\mu = [3, 2]$ it follows from Lemma~\ref{lem:char2_eo32_genus4}, and for $\mu > [3, 2]$ by Remark~\ref{rem:priessmootheo}, and for $\mu = [3, 2, 1]$ by \cite[Proposition~3.1]{stohrvoloch}.

Fix $\mu > [4, 3]$ and note that $\bZ_{\mu}$ is irreducible in $\cA_4$ by \cite[Theorem 11.5]{ekedahlvdgeer}. Using that $\cA_4$ is a smooth stack and $\cJ_4$ is a divisor in $\cA_4$, the dimension of any component $\Gamma$ of $\bZ_{\mu}\cap \cJ_4$ is either $\dim \Gamma = 9 - \sum \mu_i$ or $\dim \Gamma =  10 - \sum \mu_i$. If the latter is true for some~$\Gamma$, then $\Gamma = \bZ_{\mu} \subset \cJ_4$. This is impossible since $Z_{[4, 3]} \subseteq \bZ_{\mu}$ would result in a $3$-dimensional family of Jacobians $\cJ_C$ in $\cA_4$ with $\mu(\cJ_C) = [4, 3]$. Indeed, $Z_{[4, 3]}(\cM_4) = \o$, while the dimension of the locus of Jacobians of singular curves with Ekedahl-Oort type $[4, 3]$ is at most $2$.
\end{proof}

By \cite[Proposition 3.1]{stohrvoloch}, there are no smooth curves of genus $g = 3$ with Ekedahl-Oort type $[3, 2]$ in characteristic $2$. However, for any $g\geq 4$, we find the following result.

\begin{cor}
In characteristic $p = 2$, there is a smooth genus-$g$ curve $C$ with $\mu(C) =~[3, 2]$, for any $g \geq 4$.
\end{cor}
\begin{proof}
By Proposition \ref{prop:eo_conclusion_char2}, $Z_{[3, 2]}(\cM_4)$ is non-empty and is pure of the expected codimension~$5$. Therefore, \cite[Theorem 6.4]{pries_current_results} implies the result as explained in Remark \ref{rem:priessmootheo}.
\end{proof}

\subsection{Ekedahl-Oort stratification of $\cM_g$ in characteristic $3$}
In the remainder of this section, we examine the Ekedahl-Oort stratification of moduli spaces of smooth genus-$g$ curves for $g \geq 4$ in characteristic $p = 3$, building on insights from hyperelliptic curves.

Any hyperelliptic curve $C$ of genus $g \geq 2$ over an algebraically closed field $k$ with characteristic $p \neq 2$ can be written in the following normal form $$y^2 = f(x) = x^{2g + 1} + a_{2g}'x^{2g} + \ldots + a_2'x^2 + x,$$ with $a_i' \in k$, $2\leq i\leq 2g$, obtained by choosing $0$ and $\infty$ to be the branch points and by scaling; set $a_1' = a_{2g + 1}' = 1$ and $a_i' = 0$ for $i \not \in \{1, 2, \ldots, 2g - 1\}$. Furthermore, we denote $$(f(x))^{(p - 1)/2} = \sum_{i \geq 0}a_ix^i.$$ 
Then, using \cite[page~54]{stohrvoloch}, the Hasse-Witt matrix of $C$, that is the matrix of the Frobenius operator $F_C$, can be computed as: 
\begin{equation}
\label{eqn:hasse_witt_he_char3}
H = \begin{pmatrix}
a_{p - 1} & a_{2p - 1} &  \ldots & a_{pg - 1} \\ 
a_{p - 2} & a_{2p - 2} &  \ldots & a_{pg - 2}\\ 
\vdots & \vdots &\ddots & \vdots\\
a_{p - g} & a_{2p - g} &  \ldots & a_{pg - g} 
\end{pmatrix}.    
\end{equation}

Zhou studied the Ekedahl-Oort loci of hyperelliptic genus-$4$ curves for $p = 3$ in \cite{zhou_genus4}. Alternatively, one can directly verify that there exists a unique hyperelliptic curve with Ekedahl-Oort type $[4, 3]$ when $p = 3$. Let $k$ be an algebraically closed field with $\mathrm{char}(k) = 3$. Here, $a_i = a_i'$ holds in the previous description, which we apply in the following example.

\begin{exmp}
\label{exmp:he43}
By our description above, any hyperelliptic genus-$4$ curve $C$ over an algebraically closed field $k$ with characteristic $p = 3$ can be written in the following normal form $y^2 = f(x) = x^9 + a_8x^8 + \ldots + a_2x^2 + x,$ with $a_i \in k$ for $2\leq i\leq 8$, and one can compute the Hasse-Witt matrix $H$ of $C$ using the formula \eqref{eqn:hasse_witt_he_char3} as $$H = \begin{pmatrix}
a_2 & a_5 & a_8 &0 \\ 
1 & a_4 & a_7 & 0\\ 
0 & a_3 & a_6 & 1\\ 
0 & a_2 & a_5 & a_8
\end{pmatrix}.$$  By \eqref{eqn:43_eqcond}, $\mu(C) = [4, 3]$ is equivalent to the conditions  $\rank_k (H) = 2$ and $HH^{\sigma} = 0$, where~$H^{\sigma}$ is the matrix obtained by raising each of the entries of $H$ to the third power. The first condition is equivalent to $a_2 = a_5 = a_8 = 0$ using that the discriminant of $f$ should be non-zero, while the second one gives us $a_3 = a_4 = a_6 = a_7 = 0$. Hence, there is a unique hyperelliptic curve $$C: y^2 = x^9 + x$$ with $\mu(C) = [4, 3]$ up to isomorphism. In particular, $Z_{[4, 3]}(\cH_4)$ is a non-empty and finite set for $p = 3$.
\end{exmp}

As an application of Theorem \ref{thm:he_implies_eo}, we get the following result. Alternatively, the non-emptiness of $Z_{\mu}(\cM_4)$ for $ \mu > [4, 3] $ in characteristic $3$ follows from \cite[Theorem~1.1]{zhou_genus4}.

\begin{prop}
\label{prop:eo_conclusion_char3}
For any Young diagram $\mu$ such that $\mu > [4, 3]$, the locus $ Z_{\mu}(\mathcal{M}_4) $ is non-empty and pure of codimension $ {cd}(\mu) $ in $ \mathcal{M}_4 $ in characteristic $ 3 $.
\end{prop}

\begin{proof}
By Example \ref{exmp:he43}, $Z_{[4, 3]}(\cH_4)$ is finite. Theorem \ref{thm:he_implies_eo} implies the result. 
\end{proof}

\begin{rem}
Using Theorem \ref{thm:he_implies_eo}, Remark \ref{rem:he43finite_implies_ss3dim}, and Example \ref{exmp:he43}, we also conclude that there is a component $W$ of the supersingular locus of $\cM_4$ in characteristic $3$ with $\dim W = 3$. See also \cite[Corollary 4.4]{harashita} and \cite[Theorem 3.4]{dd2}.
\end{rem}

Now, we combine the conclusions for $p = 2$ and $p = 3$ and describe the generic behavior of the components of $V_{g - 4}(\cMbar_{g})$ for any $g \geq 4$ as an immediate consequence of  Theorem \ref{thm:main1}.
\begin{cor}
In characteristic $p \in \{ 2, 3\}$, every component of $V_{g - 4}(\cMbar_{g})$ has a generic {$a$-number} $1$, for any $g\geq 4$. 
\end{cor}
\begin{proof}
For $p = 2$ and $p = 3$, the result follows from Propositions \ref{prop:eo_conclusion_char2} and \ref{prop:eo_conclusion_char3}, respectively, using $\mu = [4, 1]$ and Theorem \ref{thm:main1} (or Theorem \ref{thm:fvgdp} (2)).
\end{proof}

\begin{rem}
This article aimed to develop techniques applicable to various Young diagrams~$\mu$, which we applied to the case $\mu = [g - 4, 1]$ in the preceding corollary. Here, we provide a more elementary proof of the corollary for $p = 3$, avoiding the complex machinery from earlier sections. By Theorem \ref{thm:main1}, it suffices to show that every component of $V_0(\cMbar_4)$ has a generic $a$-number of $1$ in characteristic $p = 3$. This remaining claim follows from two facts: every component of $V_0(\cMbar_4)$ contains a point corresponding to a smooth hyperelliptic curve (see \cite[Proposition 2.7]{fabervdgeer}), and every component of $V_0(\mathcal{H}_4)$ has a generic $a$-number $1$ in characteristic $p = 3$ (see \cite[Theorem 1.1]{zhou_genus4}).
\end{rem}

When $g = 4$, Theorem \ref{thm:he_implies_eo} indicated the benefits of studying hyperelliptic curves $C$ such that $F_C^2 = 0.$ However, for $g \geq 5$ in characteristic $3$, we show below that no hyperelliptic curves $C$ satisfy this property. Some interesting Ekedahl-Oort types are characterized by the condition $F_C^2 = 0$. Namely, $$\mu = [g, g - 1, \ldots, \lfloor g/2 \rfloor],$$ is equivalent to $\rank_k(F_C) = \lfloor g/2 \rfloor$ and $\rank_k(F_C^2) = 0$, while $$\mu = [g, \mu_2, \ldots, \mu_{g - 1}]$$  is equivalent to $\rank_k(F_C) = 1$ and $\rank_k(F_C^2) = 0$. Curves with the first Ekedahl-Oort type mentioned above belong to the largest Ekedahl-Oort stratum entirely contained in the supersingular locus, as shown in \cite[Theorem~4.8, Step 2]{chaioort}, while curves with the latter have $p$-rank $0$ and $a$-number $g - 1$.
\begin{prop}
If $C$ is a smooth hyperelliptic curve of genus $g\geq 5$ in characteristic~$3$, then $F_C^2 \neq 0$ on $H^1(C, \cO_C)$, where $F_C$ is the Frobenius operator of $C$. In particular, $Z_{\mu}(\cH_g) = \emptyset$ for at least $\mathrm{Fibonacci}(g+1)$ Young diagrams $\mu$ in characteristic $3$. 
\end{prop}

\begin{proof}
Let $g \geq 5$. We compute the the Hasse-Witt matrix $H$ of a hyperelliptic curve $C: y^2 = x^{2g + 1} + a_{2g}x^{2g} + \ldots + a_2x^2 + x$ in characteristic $3$ using \eqref{eqn:hasse_witt_he_char3}, and get $$H = \begin{pmatrix}
a_{2} & a_{5} & a_8 & \ldots & a_{3g - 4}& 0 \\ 
1 & a_{4} & a_7 & \ldots & a_{3g - 5}& 0\\ 
0 & a_{3} & a_6 & \ldots & a_{3g - 6}& 0 \\ 
0 & a_{2} & a_5 & \ldots & a_{3g - 7}& a_{3g - 4}\\ 
0 & 1 & a_4 & \ldots & a_{3g - 8}& a_{3g - 5}\\ 
\vdots & \vdots & \vdots & \ddots & \vdots & \vdots\\
0 & a_{5 - g} & a_{8 - g}&  \ldots & a_{2g - 2} & a_{2g + 1}\\
0 & a_{6 - g} & a_{9 - g}& \ldots & a_{2g - 3} & a_{2g}
\end{pmatrix}.
$$
Therefore, the $(5, 1)$-entry of $H H^{\sigma}$ equals 1, and $H H^{\sigma} \neq 0$, which implies the result. The count $\mathrm{Fibonacci}(g + 1) = \sum_{i = 0}^{\lfloor \frac{g}{2} \rfloor} {g - i \choose i}$ follows by counting the number of possible final types~$\nu$ as in \eqref{eqn:nu_def} that $C$ could have, using the property $\nu(i) \leq \nu(i + 1) \leq \nu(i) + 1$ for $i \geq 0$.
\end{proof}


\begin{thebibliography}{vdGvdV95}


\bibitem[Ale05]{alexeev}
{Alexeev, V.}, \emph{Compactified Jacobians and Torelli map.},
Publ. Res. Inst. Math. Sci., 40(4), pp.~1241–1265, 2004.

\bibitem[AP08]{achterpries}
{Achter, J. and Pries, R.}, \emph{Monodromy of the $p$-rank strata of the moduli space of curves.},
Int. Math. Res. Not. IMRN, 15, Art. ID rnn053, 25 pp., 2008.

\bibitem[AP11]{achterprieshe}
{Achter, J. and Pries, R.}, \emph{The {$p$}-rank strata of the moduli space of hyperelliptic curves.},
Adv. Math., 227, no. 5, pp.~1846-1872, 2011.


\bibitem[Bou98]{bouw}
{Bouw, I.}, \emph{The p-rank of curves and covers of curves.},
Progr. Math 187, pp. ~403-412, 1998.

\bibitem[CO11]{chaioort}
{Chai, C. and Oort, F.}, \emph{Monodromy and irreducibility of leaves.},
Ann. of Math., vol. 175, pp.~1359--1396, 2011.


\bibitem[Dia87]{diaz}
{Diaz, S.}, \emph{Complete subvarieties of the moduli space of smooth curves. In Algebraic
geometry}, In Algebraic
geometry, Bowdoin 1985 (S. Bloch, ed.), Proc. Sympos. Pure Math., 46, Part 1,
Amer. Math. Soc., Providence, RI, pp.~77–81, 1987.


\bibitem[Dra24]{dd2}
{Dragutinovi\'c, D.}, \emph{Supersingular curves of genus four in characteristic two.} Proc. Amer. Math. Soc., 152, 6, pp.~2333–2347, 2024.


\bibitem[EvdG09]{ekedahlvdgeer}
{Ekedahl, T. and van der Geer, G.}, \emph{Cycle Classes of the E-O Stratification on the Moduli of Abelian Varieties.}, In: Algebra, Arithmetic, and Geometry: Volume I: In Honor of Yu. I. Manin, pp. ~567-636. Birkhäuser Boston, Boston, 2009.



\bibitem[EP07]{elkinpries}
{Elkin, A. and Pries, R.}, \emph{Hyperelliptic curves with $a$-number 1 in small characteristic.}, Albanian Journal of Math, 1 (4), pp. ~245-252, 2007.



\bibitem[EP13]{elkinpries_he}
{Elkin, A. and Pries, R.}, \emph{Ekedahl–Oort strata of hyperelliptic curves in characteristic 2}, Algebra Number Theory, vol. 7, no. 3, pp.~507-532, 2013.


\bibitem[FvdG04]{fabervdgeer}
{Faber, C. and van der Geer, G.}, \emph{Complete subvarieties of moduli spaces and the Prym map},
 J. Reine Angew. Math., 573, pp.~117-137, 2004.

\bibitem[GP05]{glasspries}
{Glass, D. and Pries, R.}, \emph{Hyperelliptic curves with prescribed {$p$}-torsion},
Manuscripta Math., 117, no. 3, pp.~299-317, 2005.



\bibitem[Har22]{harashita}
{Harashita, S.}, \emph{Supersingular abelian varieties and curves, and their moduli spaces, with a remark on the dimension of the moduli of supersingular curves of genus 4}, Theory and Applications of Supersingular Curves and Supersingular Abelian Varieties, RIMS Kôkyûroku Bessatsu, Res. Inst. Math. Sci. (RIMS), vol. B90, pp.~1-16, 2022.


\bibitem[Knu83]{knudsen}
{Knudsen, F.}, \emph{The projectivity of the moduli space of stable curves. II. The stacks $M_{g,n}$}, Math. Scand. 52, no. 2, pp.~161–199., 1983.

\bibitem[Kob75]{koblitz}
{Koblitz, N.}, \emph{$p$-adic variation of the zeta-function over families of varieties defined over
finite fields.}, Compositio Math. 31, pp.~119-218., 1975.


\bibitem[Loo95]{looijenga}
{Looijenga, E.}, \emph{On the tautological ring of $\cM_g$.},
Invent. Math. 121, pp.~411–419., 1995.

\bibitem[Moo22]{moonen}
{Moonen, B.}, \emph{Computing discrete invariants of varieties in positive characteristic I. Ekedahl-Oort types of curves.},
J. Pure Appl. Algebra, Vol. 226, 11, Paper No. 107100, 19 pp.,  2022.

\bibitem[Oda69]{oda}
{Oda, T.}, \emph{The first de Rham cohomology group and Dieudonné modules},
Annales scientifiques de l'École Normale Supérieure, Serie 4, Volume 2, no. 1, pp. 63-135., 1969.

\bibitem[Oor62]{oort:genjac}
Oort, F., \emph{A construction of generalized Jacobian varieties by group extensions}, Mathematische Annalen,
vol. 147, pp.~277-286., 1962.

\bibitem[Oor74]{oort_subv}
Oort, F., \emph{Subvarieties of moduli spaces}, Invent. Math., 24, pp.~95-119., 1974.

\bibitem[Oor91]{oort_hess}
Oort, F., \emph{Hyperelliptic supersingular curves},
In: van der Geer, G., Oort, F., Steenbrink, J. (eds) Arithmetic Algebraic Geometry. Progress in Mathematics, vol 89. Birkhäuser, Boston, MA., pp.~247-284., 1991.

\bibitem[Oor01]{oort}
Oort, F., \emph{A stratification of a moduli space of abelian varieties},
Moduli of abelian varieties
(Texel Island, 1999), vol. 195 of Progr. Math., Birkh\"auser, Basel, pp.~345-416., 2001.

\bibitem[PWZ11]{pinkwedhornziegler}
{Pink, R. and Wedhorn, T. and Ziegler P.}, \emph{Algebraic zip data},  Doc. Math., vol. 16, pp.~253-300, 2011.


\bibitem[Pri08]{pries_eo}
Pries, R., \emph{A short guide to $p$-torsion of abelian varieties in characteristic~$p$},
Computational arithmetic geometry, Amer. Math. Soc., Providence, RI, vol. 463, pp.~121-129, 2008.

\bibitem[Pri09]{pries_a_number}
Pries, R., \emph{The $p$-torsion of curves with large $p$-rank},
Int. J. Number Theory 5, 6: pp.~1103-1115, 2009.  


\bibitem[Pri19]{pries_current_results}
Pries, R., \emph{Current results on Newton polygons of curves},
Open Problems in Arithmetic Algebraic Geometry, Oort (ed), Advanced Lectures in Mathematics, 46, ch. 6, pp.~179-208, 2019.


\bibitem[SV87]{stohrvoloch}
{Stöhr, K.-O. and Voloch, J.F.}, \emph{A formula for the \text{Cartier} operator on plane algebraic curves},
 J. Reine Angew. Math. 377, pp.~49--64, 1987.

\bibitem[vdG99]{vdgeercycle}
{van der Geer, G.}, \emph{Cycles on the moduli space of abelian varieties},
 Aspects Math., E33. {pp.~65-89}, 1999.

\bibitem[Wed05]{wedhorn} Wedhorn, T., \emph{Specialization of F-zips.}, 	arXiv:math/0507175, 2005.



\bibitem[Xar20]{xarles}
Xarles, X., {\emph{A census of all genus 4 curves over the field with 2 elements}, arXiv:2007.07822}, 
2020. Repository: \url{https://github.com/XavierXarles/Censusforgenus4curvesoverF2} 


\bibitem[Zho20]{zhou_genus4}
Zhou, Z., \emph{Ekedahl-Oort strata on the moduli space of curves of genus four.},
Rocky Mountain J. Math. 50, no. 2, pp. 747–761, 2020.




\end{thebibliography}
\end{document}